\documentclass{article}

\usepackage{amsmath,amssymb}
\usepackage[numbers]{natbib}
\usepackage{rotating}
\usepackage{graphicx}

\newcommand{\R}{\mathbb{R}}

\begin{document}

\title{Singular Quasilinear Elliptic Problems \\ With Convection Terms}
\author{
	\bf Umberto Guarnotta\\
	\small{Dipartimento di Matematica e Informatica, Universit\`a di Catania,}\\
	\small{Viale A. Doria 6, 95125 Catania, Italy}\\
	\small{\it E-mail: umberto.guarnotta@phd.unict.it}
}
\date{}
\maketitle

\begin{abstract}
In this paper we present some very recent results regarding existence, uniqueness, and multiplicity of solutions for quasilinear elliptic equations and systems, exhibiting both singular and convective reaction terms. The importance of boundary conditions (Dirichlet, Neumann, or Robin) is also discussed. Existence is achieved via sub-supersolution and truncation techniques, fixed point theory, nonlinear regularity, and set-valued analysis, while uniqueness and multiplicity are obtained by monotonicity arguments.
\end{abstract}

\section{Introduction}

Let $ \Omega \subseteq \R^N $, $ N \geq 3 $, be a bounded domain with $ C^2 $-boundary $ \partial \Omega $, and let $ n = n(x) $ the outer unit normal vector to $ \partial \Omega $ at its point $ x $. The families of problems we are dealing with can be written in the form
\begin{equation}
\label{scalar}
\left\{
\begin{array}{ll}
-{\rm div} \, a(\nabla u) + \lambda u^{p-1} = h(x,u,\nabla u) \quad &\mbox{in} \;\; \Omega, \\
u>0 \quad &\mbox{in} \;\; \Omega, \\
a(\nabla u) \cdot n + \beta u^{p-1} = 0 \quad &\mbox{on} \;\; \partial \Omega,
\end{array}
\right.
\end{equation}
and
\begin{equation}
\label{vector}
\left\{
\begin{array}{ll}
-\Delta_p u = f(x,u,v,\nabla u,\nabla v) \quad &\mbox{in} \;\; \Omega, \\
-\Delta_q v = g(x,u,v,\nabla u,\nabla v) \quad &\mbox{in} \;\; \Omega, \\
u,v>0 \quad &\mbox{in} \;\; \Omega, \\
\nabla u \cdot n = \nabla v \cdot n = 0 \quad &\mbox{on} \;\; \partial \Omega,
\end{array}
\right.
\end{equation}
concerning a single equation or a system, respectively. Here, $ 1 < p < +\infty$ and $ \lambda, \beta $ denote non-negative constants satisfying $ \lambda+\beta > 0 $; for any $ 1 < r < +\infty $, $ \Delta_r $ stands for the $ r $-Laplacian operator, that is, $ \Delta_r u := {\rm div} \, (|\nabla u|^{r-2} \nabla u) $, while $ a:\R^N \to \R^N $ is a continuous, strictly monotone map having suitable properties, which basically stem from Lieberman's nonlinear regularity theory \cite{L} and Pucci-Serrin's maximum principle \cite{PS}. \\
In order to give the idea, the operator $ u \mapsto {\rm div} \, a(\nabla u) $ is patterned after the $ (p,q) $-Laplacian, $ 1 < q < p < +\infty$, a non-homogeneous operator of the form $ \Delta_p+\Delta_q $, but it also encompasses the $ r $-Laplacian. Hereafter, we will suppose $ a(\xi) = a_0(|\xi|) \, \xi $ for any $ \xi \in \R^N $, being $ a_0: (0,+\infty) \to (0,+\infty) $ an opportune $ C^1 $ function, and denote by $ G: \R^N \to \R $ the map
\[
G(\xi) = \int_{0}^{|\xi|} \tau a_0(\tau) \, d\tau.
\]
We observe that $ \nabla G = a $ in $ \R^N $ (see, e.g., \cite{GMP}). In addition, $ h:\Omega \times (0,+\infty) \times \R^N \to [0,+\infty) $ and $ f,g: \Omega \times (0,+\infty)^2 \times \R^{2N} \to \R $ are assumed to be Carathéodory functions, which satisfy certain growth hypotheses and can be possibly both singular and convective; see \cite{GMM,GM} for further details. Just to give the flavor, one can consider the nonlinearities
\begin{eqnarray*}
h(x,s,\xi) &=& a(x)(s^{-\eta} + s^{p-1} + |\xi|^{p-1}), \\
f(x,s,t,\xi_1,\xi_2) &=& (\sin s) (s^{-\alpha_1} t^{\beta_1} - |\xi_1|^{\gamma_1} - |\xi_2|^{\delta_1}), \\
g(x,s,t,\xi_1,\xi_2) &=& (\cos t) (s^{\alpha_2} t^{-\beta_2} - |\xi_1|^{\gamma_2} - |\xi_2|^{\delta_2}), \\
\end{eqnarray*}
with $ \alpha_i,\beta_i,\gamma_i,\delta_i,\eta>0 $, $ i=1,2 $, $ \max\{\gamma_1,\delta_1\}<\beta_1-\alpha_1<p-1 $, $ \max\{\gamma_2,\delta_2\}<\alpha_2-\beta_2<q-1 $, and $ a \in L^\infty(\Omega) $, with $ a \geq 0 $ a.e. in $ \Omega $.

A problem in form (\ref{scalar}) has been studied in \cite{GMM}, while \cite{GM} analyzes (\ref{vector}). The Dirichlet boundary value problem related to (\ref{scalar}) is partially treated in \cite{LMZ}; moreover, \cite{CLM} represents the Dirichlet counterpart to (\ref{vector}). To the best of our knowledge, there are few other contributions in this direction: here we only mention \cite{MMZ}, as an example. \\
On the other hand, singular problems and convective ones have been intensively studied in the last decade, although separately: the reader who is interested in singular problems can consult \cite{PW} for equations, as well as \cite{MM} for systems, and the references therein; concerning convective equations, we refer to the recent papers \cite{FMM,FMP}.

\section{The technique}

There are two main difficulties to overcome: (i) the loss of variational structure, due to the presence of convection terms; (ii) the singularity of reaction terms. First of all, let us discuss (\ref{scalar}) with a particularized nonlinearity $ h $, splitted as the sum of a convective term and a singular one, that is,
\begin{equation}
\label{splitreact}
h(x,s,\xi) = f(x,s,\xi) + g(x,s),
\end{equation}
with $ f: \Omega \times \R \times \R^N \to [0,+\infty), g: \Omega \times (0,+\infty) \to [0,+\infty) $ being Carathéodory functions; it is worth noticing explicitly that $ f(x,\cdot,\xi) $ is continuous on the whole $ \R $, for a.a. $ x \in \Omega $ and all $ \xi \in \R^N $, so the convection term is non-singular. In order to recover the variational structure of the problem, we can `freeze' the gradient term in the differential equation, obtaining a one-parameter family of problems depending on $ w \in C^1(\overline{\Omega}) $:
\begin{equation}
\label{auxscalar}
\left\{
\begin{array}{ll}
-{\rm div} \, a(\nabla u) = f(x,u,\nabla w) + g(x,u) \quad &\mbox{in} \;\; \Omega, \\
u>0 \quad &\mbox{in} \;\; \Omega, \\
a(\nabla u) \cdot n + \beta |u|^{p-2}u = 0 \quad &\mbox{on} \;\; \partial \Omega,
\end{array}
\right.
\end{equation}
where we have assumed $ \lambda = 0 < \beta $ for the sake of simplicity. \\
The situation looks quite different for system (\ref{vector}), for two reasons: the equations are coupled, and the corresponding operators are not maximal monotone, because of Neumann boundary conditions and the absence of potential terms $ u^{p-1} $, $ v^{q-1} $. In this case, after adding them, the `freezing technique' has to be applied to the right-hand side of the system, instead of only gradient terms; in particular, the auxiliary problem to solve is
\begin{equation}
\label{auxvector}
\left\{
\begin{array}{ll}
-\Delta_p u + |u|^{p-2}u = f(x,z_1,z_2,\nabla w_1,\nabla w_2) + z_1^{p-1} \quad &\mbox{in} \;\; \Omega, \\
-\Delta_q v + |v|^{q-2}v = g(x,z_1,z_2,\nabla w_1,\nabla w_2) + z_2^{q-1} \quad &\mbox{in} \;\; \Omega, \\
\nabla u \cdot n = \nabla v \cdot n = 0 \quad &\mbox{on} \;\; \partial \Omega,
\end{array}
\right.
\end{equation}
where $ (z_1,z_2,w_1,w_2) \in W^{1,p}(\Omega) \times W^{1,q}(\Omega) \times C^1(\overline{\Omega})^2 $, with $ z_1,z_2 > 0 $ in $ \Omega $, is fixed.

Now, according to the boundary conditions in (\ref{auxscalar}), we are able to define a suitable norm $ \|\cdot\| $ in $ W^{1,p}(\Omega) $ such that the operator $ u \mapsto -{\rm div} \, a(\nabla u) $ admits a potential. Indeed, let
\[
\|u\| := \left( p \int_\Omega G(\nabla u) dx + \beta \int_{\partial \Omega} |u|^p d\sigma \right)^{\frac{1}{p}}
\]
for any $ u \in W^{1,p}(\Omega) $, where $ d\sigma $ is the Hausdorff $ (N-1) $-dimensional measure (and $ u_{\mid_{\partial \Omega}} $ has to be understood in the sense of traces). Robin boundary conditions and nonlinear Green's formula yield
\begin{equation*}
\begin{split}
&\left\langle D_F \left( \frac{1}{p} \|u\|^p \right), \phi \right\rangle = \int_\Omega a(\nabla u) \cdot \nabla \phi dx + \beta \int_{\partial \Omega} |u|^{p-2}u \phi d\sigma \\
&= \int_\Omega a(\nabla u) \cdot \nabla \phi dx - \int_{\partial \Omega} (a(\nabla u) \cdot n) \phi d\sigma = \langle -{\rm div} \, a(\nabla u), \phi \rangle,
\end{split}
\end{equation*}
for any $ \phi \in W^{1,p}(\Omega) $, being $ D_F $ the Fréchet derivative and $ \langle \cdot, \cdot \rangle $ the standard duality brackets. A similar argument can be used for system (\ref{auxvector}), choosing the norms
\[
\|u\|:= \left( \int_\Omega |\nabla u|^p dx + \int_{\Omega} |u|^p dx \right)^{\frac{1}{p}} \quad \mbox{and} \quad \|v\|:= \left( \int_\Omega |\nabla v|^q dx + \int_{\Omega} |v|^q dx \right)^{\frac{1}{q}}
\]
in $ W^{1,p}(\Omega) $ and $ W^{1,q}(\Omega) $ respectively, and recalling that $ \nabla z \cdot n = 0 $ if and only if $ |\nabla z|^{r-2} (\nabla z \cdot n) = 0 $ (for any $ 1 < r < +\infty $). In Dirichlet problems exhibiting no potential terms, Poincaré's inequality permits to define a suitable norm, as above.

The right-hand sides of (\ref{auxscalar})-(\ref{auxvector}) can be opportunely modified to fit the classical variational pattern: it suffices to control the singular behavior. To do this, we can employ the truncation method, provided the reaction terms grow opportunely. Firstly, we have to construct a subsolution; hence, we will define another auxiliary problem, which is variational, solve it through the direct methods of calculus of variations, and prove - via comparison arguments - that the found solution actually is a solution to (\ref{auxscalar}) (or (\ref{auxvector})). To better explain this crucial passage, we argue for (\ref{auxscalar}). The term $ f(x,u,\nabla w) $ is standard, so no problem arises; on the other hand, there are some ways to guarantee the existence of a subsolution. One of them consists in requiring $ g(x,\cdot) $ to be a singular function (see, e.g., \cite{CGP} for details), that is,
\[
\lim_{s \to 0^+} g(x,s) = +\infty \quad \mbox{uniformly w.r.t.} \;\; x \in \Omega,
\]
or requiring a monotonicity condition (cf. \cite{LMZ,GMM}), as
\[
g(x,\cdot) \;\; \mbox{is non-increasing in} \; (0,1], \quad g(\cdot,1) \not\equiv 0.
\]
An extra hypothesis is also required: for instance,

\begin{equation}
\label{growth}
g(x,s) \leq Cs^{-\gamma} \quad \forall (x,s) \in \Omega \times (0,1),
\end{equation}
being $ C > 0, \gamma \in (0,1) $ suitable constants. For Dirichlet problems (see \cite{CLM} for systems), (\ref{growth}) can be regarded as a summability hypothesis on the superposition operator $ x \mapsto g(x,u(x)) $, according to Hardy-Sobolev's inequality: for any $ \phi \in W^{1,p}_0(\Omega) $ and $ u $ satisfying $ u(x) \geq k {\rm dist}(x,\partial \Omega) $ for any $ x \in \Omega $ and some $ k > 0 $, we get
\[
\begin{split}
\int_\Omega g(x,u(x))|\phi(x)| dx &\leq C \int_\Omega u(x)^{-\gamma}|\phi(x)| dx \\
&\leq Ck^{-\gamma} \int_{\Omega} ({\rm dist}(x,\partial \Omega))^{-\gamma} |\phi(x)| dx \\
&\leq C' \int_\Omega |\nabla \phi(x)|^p dx < +\infty,
\end{split}
\]
being $ C' > 0 $ opportune. A different summability hypothesis can be found in \cite{GMM}. On the contrary, for Neumann problems (even systems) the situation looks easier, because adding a constant does not affect neither first nor second derivatives (vide \cite{GM}). \\
For the sake of completeness, we mention another technique to avoid singularities (which is used, e.g., in \cite{FMM}): the underlying idea is to `shift' the values of the singular term, and to solve the $ \epsilon $-dependent problems obtained by substituting $ g(x,u) $ with $ g(x,u+\epsilon) $; hence, trying to pass the limit in $ \epsilon $, one might find a solution to (\ref{auxscalar}).

Denoting with $ \underline{u} $ the constructed subsolution, and letting $ T:W^{1,p}(\Omega) \to W^{1,p}(\Omega) $ the truncation operator
$ T(u) = \max \left\{u,\underline{u}\right\} $, problem (\ref{auxscalar}) can be transformed into
\[
\left\{
\begin{array}{ll}
-{\rm div} \, a(\nabla u) = f(x,u,\nabla w) + g(x,T(u)) \quad &\mbox{in} \;\; \Omega, \\
u>0 \quad &\mbox{in} \;\; \Omega, \\
a(\nabla u) \cdot n + \beta |u|^{p-2}u = 0 \quad &\mbox{on} \;\; \partial \Omega.
\end{array}
\right.
\]
A solution of this new problem, say $ u $, can be found using Weierstrass-Tonelli's theorem, and a simple comparison argument reveals that $ u \geq \underline{u} $, so we have solved (\ref{auxscalar}). System (\ref{auxvector}) can be treated in a similar way.

At this point, we are ready to start the `unfreezing procedure': given the map $ \Psi(w) = u $, being $ u $ solution to (\ref{auxscalar}), we would like to find $ u^* $ such that $ \Psi(u^*) = u^* $, and hence $ u^* $ turns out to be a solution to (\ref{scalar}). This is a fixed point problem, but unfortunately we are not able to guarantee that $ \Psi $ satisfies the hypotheses of any fixed point theorem. Hence, we consider the multi-function $ \mathcal{S}: C^1(\overline{\Omega}) \to 2^{C^1(\overline{\Omega})} $ defined by 
\[
\mathcal{S}(w) := \left\{u \in C^1(\overline{\Omega}): \, u \; \mbox{is a solution to} \; (\ref{auxscalar}), \; u \geq \underline{u} \right\},
\]
and we show that $ S(w) $ is downward directed for all $ w \in C^1(\overline{\Omega}) $, in order to well-define its selection $ \mathcal{T}(w) := \min \mathcal{S}(w) $. Now compactness and continuity of $ \mathcal{T} $ are basically inherited by compactness and lower semicontinuity of $ \mathcal{S} $ (the proof of lower semicontinuity is rather technical, and based on an approximation procedure: we refer to \cite{FMP,LMZ,GMM}). An application of Schaefer's fixed point theorem (see \cite[p. 827]{GP}) on $ \mathcal{T} $, together with suitable algebraic conditions on the growth parameters of $ f $ and $ g $, ensures that there exists $ u^* $ such that $ \mathcal{T}(u^*) = u^* $, and so we are done. \\
The procedure is more delicate when we discuss about systems: we have both $ (z_1,z_2) $ and $ (w_1,w_2) $ to `unfreeze', but the two problems can be handled with a fixed point approach again. Constructing a supersolution in order to gain compactness, besides modifying the truncated problem, allows to apply Schauder's fixed point theorem to get a new problem depending only on $ (w_1,w_2) $. Unfortunately, we cannot apply Schaefer's theorem: as the operator is not maximal monotone, we are not able to prove the boundedness of the set of solutions to the problem: $ t \mathcal{T}(u) = u $ for some $ t \in (0,1) $. Hence, we need the following a priori estimates on the gradients (cf. \cite[Theorem 3.1]{CiM}):
\[
\begin{split}
\|\nabla u\|_{L^\infty(\Omega)} &\leq C \|f(\cdot,u,v,\nabla w_1,\nabla w_2)\|_{L^\infty(\Omega)}^{\frac{1}{p-1}}, \\
\|\nabla v\|_{L^\infty(\Omega)} &\leq C \|g(\cdot,u,v,\nabla w_1,\nabla w_2)\|_{L^\infty(\Omega)}^{\frac{1}{q-1}}.
\end{split}
\]
Now we take $ D \subseteq C^1(\overline{\Omega})^2 $ such that $ (\underline{u},\underline{v}) \leq (u,v) \leq (\overline{u},\overline{v}) $ and $$ \max \{\|\nabla u\|_{L^\infty(\Omega)},\|\nabla v\|_{L^\infty(\Omega)}\} \leq M $$ for any $ (u,v) \in D $. For a large $ M > 0 $, the restriction $ \mathcal{S}_{\mid_D} $ maps $ D $ into $ 2^D $; moreover, boundedness of both gradients in $ D $ and sub-super-solutions permit to apply Schauder's theorem once again, concluding the proof. This method is related to the `trapping region' argument: cf. \cite{CaM}. Incidentally, working with a trapping region allows to consider supercritical nonlinearities or reaction terms which `roughly' blow up, in the sense that they do not satisfy (\ref{growth}); see \cite{GM} for some examples.

Concluding, it is worth spending few words about uniqueness and multiplicity of the solutions found above. Usually, linear problems possess only one solution, whereas multiplicity is encountered in nonlinear phenomena. When $ p = 2 $, problem (\ref{scalar}) admits a unique solution, under some Lipschitz-type conditions on nonlinearities $ f $ and $ g $ (cf. \cite{GMM}). The problem seems to be open for $ p \neq 2 $. \\
On the other hand, Neumann problems naturally possess multiple solutions: it suffices to think about the Laplace equation $ -\Delta u = f(x) $ under Neumann boundary condition, which possesses the family of solutions $ \{u+c: \, c \in \R\} $. However, this is not so obvious when the reaction term depends on the solution. In this spirit, we construct a (pointwise) ordered sequence of sub-super-solutions $ \underline{u}_1 < \overline{u}_1 < \underline{u}_2 < \overline{u}_2 < \ldots $, and we find infinitely many solutions $ \{u_n\}_n $ to (\ref{vector}) such that $ \underline{u}_n \leq u_n \leq \overline{u}_n $ for any $ n $ (see \cite{GM}).

\bibliographystyle{aipnum-cp}

\end{document}